\documentclass[11pt]{amsart}
\usepackage{amsmath, amsfonts, upref}
\newtheorem{prop}{Proposition}[section]

\newtheorem{lemma}[prop]{Lemma}

\newtheorem{theorem}[prop]{Theorem}

\begin{document}
\newcommand{\ee}{\mathbb{E}}
\newcommand{\pp}{\mathbb{P}}
\newcommand{\ii}{\mathbb{I}}
\newcommand{\zz}{\mathbb{Z}}
\newcommand{\rr}{\mathbb{R}}
\newcommand{\ex}{\mathrm{Exp}(1)}
\newcommand{\ra}{\rightarrow}
\newcommand{\fp}{f^\prime}
\newcommand{\fpp}{f^{\prime\prime}}
\newcommand{\ww}{W^\prime}
\newcommand{\frc}[2]{{\textstyle{\frac{#1}{#2}}}}
\newtheorem{thm}{Theorem}
\newtheorem{lmm}{Lemma}

\begin{center}
{\bf An Inductive Proof of the Berry-Esseen Theorem for Character Ratios}
\end{center}

\begin{center}
{\bf Running head: Berry-Esseen Theorem for Character Ratios}
\end{center}

\begin{center}
Submitted 3/9/05; Revised 8/6/06 
\end{center}

\begin{center}
By Jason Fulman
\end{center}
\begin{center}
Department of Mathematics, University of Southern California
\end{center}
\begin{center}
Los Angeles, CA 90089, USA
\end{center}

\begin{center}
fulman@usc.edu
\end{center}

{\bf Abstract}: Bolthausen used a variation of Stein's method to give
an inductive proof of the Berry-Esseen theorem for sums of independent,
identically distributed random variables. We modify this technique to
prove a Berry-Esseen theorem for character ratios of a random
representation of the symmetric group on transpositions. An analogous
result is proved for Jack measure on partitions.

\begin{center}
\end{center}

\begin{center}
\end{center}

\begin{center}
\end{center}

\section{Introduction}

	The Plancherel measure of a finite group $G$ is a probability
measure on the set of irreducible representations of $G$ which chooses
a representation $\rho$ with probability $\frac{dim(\rho)^2}{|G|}$,
where $dim(\rho)$ denotes the dimension of $\rho$. For instance if $G$
is the symmetric group, the irreducible representations are
parameterized by partitions $\lambda$ of $n$, and the Plancherel
measure chooses a partition $\lambda$ with probability
$\frac{n!}{\prod_{x \in \lambda} h(x)^2}$ where the product is over
boxes in the partition and $h(x)$ is the hooklength of a box. The
hooklength of a box $x$ is defined as 1 + number of boxes in same row
as x and to the right of x + number of boxes in same column of x and
below x. For example we have filled in each box in the partition of 7
below with its hooklength \[ \begin{array}{c c c c} \framebox{6}&
\framebox{4}& \framebox{2}& \framebox{1}\\ \framebox{3}&
\framebox{1}&& \\ \framebox{1} &&& \end{array} \] and the Plancherel
measure would choose this partition with probability
$\frac{7!}{(6*4*3*2)^2}$. Recently there has been interest in the
statistical properties of partitions chosen from Plancherel measure
and we refer the reader to the surveys \cite{AD}, \cite{De} and the
seminal papers \cite{J}, \cite{O1}, \cite{BOO} for a glimpse of the
remarkable recent work on Plancherel measure. We recommend \cite{Sa}
as an introduction to representation theory of the symmetric group.

	Let $\lambda$ be a partition of $n$ chosen from the Plancherel
        measure of the symmetric group $S_n$ and let
        $\chi^{\lambda}(12)$ be the irreducible character
        parameterized by $\lambda$ evaluated on the transposition
        $(12)$. The quantity $\frac{\chi^{\lambda}(12)}{dim(\lambda)}$
        is called a character ratio and is crucial for analyzing the
        convergence rate of the random walk on the symmetric group
        generated by transpositions \cite{DS}. In fact Diaconis and
        Shahshahani prove that the eigenvalues for this random walk
        are the character ratios
        $\frac{\chi^{\lambda}(12)}{dim(\lambda)}$ each occurring with
        multiplicity $dim(\lambda)^2$. Character ratios on
        transpositions also play an essential role in work on the
        moduli space of curves \cite{EO}, \cite{OP}.

	Given these motivations, it is natural to study the
 distribution of the character ratio
 $\frac{\chi^{\lambda}(12)}{dim(\lambda)}$ and there has been a
 substantial amount of work in this direction, which we now
 summarize. Kerov \cite{K1} proved that if $\lambda$ is chosen from
 the Plancherel measure of the symmetric group, then for all real
 $x_0$, \[ lim_{n \rightarrow \infty} \pp \left( \frac{n-1}{\sqrt{2}}
 \frac{\chi^{\lambda}(12)}{dim(\lambda)} \leq x_0 \right) =
 \frac{1}{\sqrt{2 \pi}} \int_{-\infty}^{x_0} e^{-\frac{t^2}{2}} dt
  .\] The details of Kerov's argument appeared in \cite{IO},
 which gave a beautiful development of Kerov's work. Hora \cite{Ho}
 gave another proof of Kerov's result, exploiting the fact that the
 kth moment of a Plancherel distributed character ratio is equal to
 the chance that the random walk generated by random transpositions is
 at the identity after k steps. Both of these proofs were essentially
 combinatorial in nature and used the method of moments (and so
 information about all moments of the character ratio). Recent work of
 Sniady \cite{Sn1}, \cite{Sn2} understands these moments in terms of
 the genus expansion from random matrix theory.

	A more probabilistic approach to Kerov's result appeared in
	\cite{F1}, which proved that for all $n \geq 2$ and real
	$x_0$, \[ \left| \pp \left( \frac{n-1}{\sqrt{2}}
	\frac{\chi^{\lambda}(12)}{dim(\lambda)} \leq x_0 \right) -
	\frac{1}{\sqrt{2 \pi}} \int_{-\infty}^{x_0} e^{-
	\frac{t^2}{2}} dt \right| \leq 40.1 n^{-1/4}. \] The proof
	used Stein's method (which is fundamentally different from the
	method of moments as it only uses information about a few
	lower order moments) and random walk on the set of irreducible
	representations of the symmetric group. Note that unlike
	Kerov's original result, this result includes an error
	term. The paper \cite{F3} used martingale theory to sharpen
	the error term to $C_s n^{-s}$ for any $s < \frac{1}{2}$ where
	$C_s$ is a constant depending on $s$. The paper \cite{SS}
	developed a refinement of Stein's method which led to a proof
	of the conjecture of \cite{F1} that an error term of $C
	n^{-1/2}$ holds where $C$ is a universal
	constant.

	The purpose of the present paper is to use a completely
	different technique to prove the $C n^{-1/2}$ bound. The
	method is based on Bolthausen's \cite{Bol} ingenious inductive
	proof of the Berry-Esseen theorem for sums of independent
	identically distributed random variables. As in \cite{F3}, we
	write the character ratio as a sum of martingale differences,
	but these are neither independent nor identically distributed
	so some subtle combinatorics is required to adapt Bolthausen's
	method. This is not the first example of adapting Bolthausen's
	method to the non i.i.d. case; Bolthausen \cite{Bol} used the
	approach to study the distribution of $\sum_{1 \leq i \leq n}
	A_{i \pi(i)}$ where $A$ is a fixed $n \times n$ matrix and
	$\pi$ is a random permutation on $n$ symbols. But the case of
	character ratios is of considerable interest and quite unlike
	any other example to which his technique has been applied.

	Note that using the method of moments, a central limit theorem
is known for character ratios on the conjugacy class of i-cycles in the
symmetric group, where i is fixed \cite{K1}, \cite{IO}, \cite{Ho}. The
preprint \cite{F4}, written after this paper, uses the ``exchangeable
pairs'' version of Stein's method to obtain an $O(n^{-1/4})$ error
term for the class of i-cycles. It also gives analogs for other
algebraic structures: Gelfand pairs, twisted Gelfand pairs, and
association schemes. It would be interesting to extend the technique
of this paper to the case of i-cycles. Significant work would be
involved in doing this, since the proof of the central limit theorem
uses the fact that in Section \ref{plancherel} when we write $\frac{{n
\choose 2} \chi^{\lambda}(12)}{dim(\lambda)}$ as a sum of martingale
differences, the expected value of the square of a summand given the
previous summands is constant. This is false for general conjugacy
classes. Also it is a nontrivial combinatorial problem to give upper
bounds on the expected absolute value of the cubes of the
summands. Fortunately for the case of transpositions this can be done
without much difficulty. And the case of transpositions does seem to
have unique practical importance \cite{EO}, \cite{OP}.

	The contents of this paper are as follows. Section
	\ref{plancherel} develops the combinatorics needed to adapt
	Bolthausen's method to the case of character ratios, and then
	proves an upper bound of $C n^{-1/2}$. Section \ref{jack} then
	recalls the Jack$_{\alpha}$ measure on partitions (here
	$\alpha > 0$ is a parameter) and why it is interesting. It
	then briefly indicates the modifications to the Plancherel
	case needed to prove a central limit theorem with an error
	term of $C_{\alpha} n^{-1/2}$, where $C_{\alpha}$ is a
	constant depending on $\alpha$. This organization is natural
	since many algebraically inclined readers will want to
	understand the result for character ratios without needing
	combinatorics of Jack polynomials; thus a useful lemma is
	given an algebraic proof in Section \ref{plancherel} and a
	combinatorial proof in Section \ref{jack}.

\section{Central limit theorem for Plancherel measure} \label{plancherel}

	The random variable we wish to study is $T_n(\lambda) =
	\frac{\sqrt{n \choose 2} \chi^{\lambda}(12)}{dim(\lambda)}$
	where $\lambda$ is chosen from the Plancherel measure of the
	symmetric group $S_n$. To begin we write $T_n$ as a sum of
	other random variables. For this we need Kerov's growth
	process on partitions \cite{K2}; this has a natural
	generalization to arbitrary finite groups \cite{F3}, but we
	only recall it in the case of interest. Given a partition
	$\lambda(j)$ of size $j$, one obtains a partition
	$\lambda(j+1)$ of size $j+1$ by choosing $\lambda(j+1)$ with
	probability $\frac{dim(\lambda(j+1))}{(j+1) dim(\lambda(j))}$
	if $\lambda(j+1)$ can be obtained from $\lambda(j)$ by adding
	a single box, and with probability 0 otherwise. Thus starting
	from $\lambda(1)$, the unique partition of size $1$, one
	obtains a random sequence $(\lambda(1),\cdots,\lambda(n))$ of
	partitions. Kerov \cite{K2} proves that each $\lambda(j)$ is
	distributed according to the Plancherel measure of $S_j$.

	Given Kerov's growth process, one can write $T_n =
	\frac{1}{\sqrt{{n \choose 2}}} (X_1+ \cdots +X_n)$ where
	$X_1=0$, $\chi^{\lambda(1)}(12)$ is defined as 0, and \[ X_j=
	\frac{{j \choose 2} \chi^{\lambda(j)}(12)}{dim(\lambda(j))} -
	\frac{{j-1 \choose 2}
	\chi^{\lambda(j-1)}(12)}{dim(\lambda(j-1))} \] for $j \geq 2$.

	Lemma \ref{ismart} states that the $X_j$ are martingale
	differences satisfying special properties. We remark that
	\cite{F3} extends this lemma to more general conjugacy classes
	and groups. The notation $\ee(A|\cdot)$ means the expected
	value of $A$ given $\cdot$.

\begin{lemma} \label{ismart} (\cite{F3})
\begin{enumerate}
\item $\ee(X_j|\lambda(j-1))=0$ for $2 \leq j \leq n$ and all partitions $\lambda(j-1)$.
\item $\ee(X_j|T_n)=\frac{j-1}{\sqrt{{n \choose 2}}} T_n$ for all $1
\leq j \leq n$.
\item $\ee(X_j^2) = j-1$.
\item $\ee(T_n^2) = 1$.
\end{enumerate}
\end{lemma}

 Frobenius \cite{Fr} found the following explicit formula for the
	character ratio of the symmetric group on transpositions: \[
	\frac{\chi^{\lambda}(12)}{dim(\lambda)} = \frac{1}{{n \choose
	2}} \sum_i \left({\lambda_i \choose 2} - {\lambda_i' \choose
	2} \right)\] where $\lambda_i$ is the length of row $i$ of
	$\lambda$ and $\lambda_i'$ is the length of column $i$ of
	$\lambda$. From his formula it follows that $X_j = c(x)$ where
	$x$ is the box added to $\lambda(j-1)$ to obtain $\lambda(j)$
	and the ``content'' $c(x)$ of a box is defined as column
	number of box - row number of box.

	Lemma \ref{moments} gives the conditional second and fourth
	moments of the $X_j$'s. We emphasize that these were not
	derived or even stated in terms of character ratios, but
	rather were proved in a completely combinatorial way by
	studying the behavior of the moments of $c(x)$ where $x$ is
	the box added during Kerov's growth process. We remark that
	for other conjugacy classes, there is not an analog of the
	fact that $\ee(X_j^2|\lambda(j-1))$ is independent of
	$\lambda(j-1)$.

\begin{lemma} \label{moments} Let $\lambda(j-1)$ be a partition
of size $j-1 \geq 1$.

\begin{enumerate}
\item (\cite{K3}) $\ee(X_j^2|\lambda(j-1)) = j-1$.

\item (\cite{La}) $\ee(X_j^4|\lambda(j-1)) = {j \choose 2} + 3 \sum_{x
\in \lambda(j-1)} c(x)^2$.
\end{enumerate}
\end{lemma}

	Lemma \ref{identity} is a useful identity. Although a
	combinatorial proof can be given using properties of Schur
	functions, we defer combinatorial arguments to the more
	general setting of Jack polynomials in Section \ref{jack} and
	give an algebraic proof.

\begin{lemma} \label{identity} Let $e_r(z_1,\cdots,z_n)= \sum_{1 \leq i_1<\cdots < i_r \leq n}
z_{i_1} \cdots z_{i_r}$ be the rth elementary symmetric function
of $z_1,\cdots,z_n$. For $\lambda$ a partition of $n$, let
$e_r(\lambda)$ denote the rth elementary symmetric function of the
contents of the boxes of $\lambda$. Then $\ee(e_r(\lambda))=0$ for $1
\leq r \leq n$. \end{lemma}

\begin{proof} If $r=n$ the result is clear since the box in the first row and column
 of $\lambda$ has content 0, so that $e_n(\lambda)=0$ for all $\lambda$.

	For $1 \leq r <n$, we use the theory of Murphy elements
 \cite{Mu}; a friendly reference giving background on these elements
 is \cite{DG}. For $2 \leq i \leq n$, the ith Murphy element is
 defined as the sum of transpositions $R_i = \sum_{1 \leq j < i}
 (j,i)$. Let $z$ be the element of the group algebra of $S_n$ which is
 the sum of all permutations with $n-r$ cycles. By Proposition 2.1 of
 \cite{DG}, $z$ is the rth elementary symmetric function of the
 elements $R_2,\cdots,R_n$.

	 Since the elements $R_2,\cdots,R_n$ are simultaneously
 diagonalizable in every irreducible representation of the symmetric
 group, it follows from Murphy's determination of their eigenvalues
 that in the representation of $S_n$ parameterized by $\lambda$, $z$
 is a scalar multiple of the $dim(\lambda) \times dim(\lambda)$
 identity matrix with scalar equal to $e_r(\lambda)$. In the regular
 representation of $S_n$ the irreducible representation parameterized
 by $\lambda$ occurs with multiplicity $dim(\lambda)$. Hence the trace
 of $z$ in the regular representation is $n! \ee(e_r(\lambda))$. But
 the coefficient of the identity in $z$ is 0, so the trace of $z$ in
 the regular representation is 0, implying the result. \end{proof}

	Lemma \ref{boundmom} gives upper bounds for $\ee(|X_n|^3)$ and
	for $\ee(|T_{n-1}||X_n|^3)$. One could prove a similar bound
	(with slightly worse constants) using the concentration
	inequality for $X_n$ in the proof of Theorem \ref{main1}.

\begin{lemma} \label{boundmom} Suppose that $n \geq 3$.
\begin{enumerate}
\item $\ee(|X_n|^3) \leq (n-1) \sqrt{2n-3}$.
\item $\ee(|T_{n-1}||X_n|^3) \leq (n-1) \sqrt{2n-3}$.
\end{enumerate}
\end{lemma}

\begin{proof} By the Cauchy-Schwarz inequality,
 $\ee(|X_n|^3) \leq \sqrt{\ee(X_n^2) \ee(X_n^4)}$. By Lemma
\ref{moments}, $\ee(X_n^2) = n-1$ and \[ \ee(X_n^4) =
\ee(\ee(X_n^4|\lambda(n-1))) = {n \choose 2} + 3 \ee \left( \sum_{x
\in \lambda(n-1)} c(x)^2 \right).\] By Lemma \ref{identity} with $r=2$
and then part 4 of Lemma \ref{ismart}, \begin{eqnarray*}
\ee \left( \sum_{x \in \lambda(n-1)} c(x)^2 \right) & = & \ee \left[ \left( \sum_{x \in
\lambda(n-1)} c(x) \right)^2 - 2 e_2(\lambda(n-1)) \right]\\ & = & {n-1 \choose
2} \ee ( T_{n-1}^2 )\\ & = & {n-1 \choose 2}. \end{eqnarray*} This proves the
first assertion.

	For the second assertion, note (using part 4 of Lemma
\ref{ismart} in the final equality) that \begin{eqnarray*}
\ee(|T_{n-1}||X_n|^3) & = &\ee(\ee(|T_{n-1}||X_n|^3|\lambda(n-1)))\\ &
= & \ee(|T_{n-1}| \ee(|X_n|^3|\lambda(n-1)))\\ & \leq &
\sqrt{\ee(T_{n-1}^2) \ee(\ee(|X_n|^3|\lambda(n-1))^2)}\\ & = &
\sqrt{ \ee(\ee(|X_n|^3|\lambda(n-1))^2)}. \end{eqnarray*} The
conditional version of the Cauchy-Schwarz inequality and part 1 of
Lemma \ref{moments} give that $\ee(|X_n|^3|\lambda(n-1))^2$ is at most
\[ \ee(X_n^2|\lambda(n-1)) \ee(X_n^4|\lambda(n-1)) = (n-1)
\ee(X_n^4|\lambda(n-1)).\] Thus
\[ \sqrt{ \ee(\ee(|X_n|^3|\lambda(n-1))^2)} \leq \sqrt{(n-1) \ee(X_n^4)},\]
and the proof of the first assertion showed this to equal $(n-1) \sqrt{2n-3}$,
as desired. \end{proof}

	Now we adapt Bolthausen's \cite{Bol} inductive proof of the
	Berry-Esseen theorem for i.i.d. random variables to the
	setting of character ratios. We remark that the unpublished
	notes of Mann \cite{Man} are a useful exposition of
	Bolthausen's proof and we refer to them in the proof of Theorem
	\ref{main1}.

\begin{theorem} \label{main1} Let
 $\lambda$ be chosen from the Plancherel measure on partitions of size
	$n$.  Then for all $n \geq 2$ and real $x_0$, \[ \left| \pp
	\left( T_n(\lambda) \leq x_0 \right) - \frac{1}{\sqrt{2 \pi}}
	\int_{-\infty}^{x_0} e^{- \frac{t^2}{2}} dt \right| \leq C
	n^{-1/2}, \] where $C$ is a universal constant.
\end{theorem}

\begin{proof} It is sufficient to prove the result for $n \geq 3$, so we assume this.

	For $z$ real, let $h_{z,0}=\ii_{(-\infty,z]}$ be the indicator
	function of the set $(-\infty,z]$. For $z$ real and $b>0$, let
	$h_{z,b}$ be the function which is 1 for $x \leq z$ and then
	drops linearly to the value 0 at $z+b$ and is 0 for $x \geq
	z+b$. Let \[ \delta(b,n) = \sup_{z} \{| \ee(h_{z,b}(T_n)) -
	\Phi h_{z,b}| \} \] where $\Phi f$ is the expected value of a
	function f under the normal distribution. Note that our
	ultimate goal is to upper bound $\delta(0,n)$.

	As in Stein's method \cite{Stn}, let \[ f(x)=f_{z,b}(x) =
	e^{x^2/2} \int_{-\infty}^x (h_{z,b}(w) - \Phi h_{z,b})
	e^{-w^2/2} dw.\] Then $f'(x) - x f(x) = h_{z,b}(x) - \Phi
	h_{z,b}$, so that \[ \ee(h_{z,b}(T_n))-\Phi h_{z,b} = \ee[
	f'(T_n) - T_n f(T_n)] .\] Part 2 of Lemma \ref{ismart} with
	$j=n$ implies that \[ \ee(X_n f(T_n)) = \ee[f(T_n)
	\ee(X_n|T_n)] = \frac{n-1}{\sqrt{{n \choose 2}}} \ee(T_n
	f(T_n)), \] so that \[  \ee[
	f'(T_n) - T_n f(T_n)] = \ee \left[ f'(T_n) -
	\frac{\sqrt{{n \choose 2}}}{n-1} X_n f(T_n) \right].\] By part
	1 of Lemma \ref{ismart} and part 1 of Lemma \ref{moments},
	this is equal to \begin{eqnarray*} & & \ee \left[ f'(T_n) \right] + \ee \left[ 	\frac{X_n^2}{n-1} f'(\sqrt{\frac{n-2}{n}} T_{n-1}) -
	f'(\sqrt{\frac{n-2}{n}} T_{n-1}) \right]  \\ & & - \ee \left[
	\frac{\sqrt{{n \choose 2}}}{n-1} X_n f(T_n) - \frac{\sqrt{{n
	\choose 2}}}{n-1} X_n f(\sqrt{\frac{n-2}{n}} T_{n-1}) 
 \right]\\ & = & \ee
	\left[ f'(T_n) - f'(\sqrt{\frac{n-2}{n}} T_{n-1}) \right]\\ & & - \ee \left[
	\frac{X_n^2}{n-1} \int_0^1  f'(\sqrt{\frac{n-2}{n}}
	T_{n-1} + t \frac{X_n}{\sqrt{{n \choose 2}}}) -
	f'(\sqrt{\frac{n-2}{n}} T_{n-1}) dt \right]. \end{eqnarray*}

	Next we upper bound $\ee[ f'(T_n) - f'(\sqrt{\frac{n-2}{n}}
	T_{n-1})]$. Recall from \cite{Bol} or \cite{Man} that for any
	$x$ and $\Delta$, \[ |f'(x+\Delta) - f'(x)| \leq |\Delta| \left( 3+
	2|x| + \frac{1}{b} \int_{0}^1 \ii_{[z,z+b]}(x+s \Delta) ds \right).\]
	Thus $ \ee[ f'(T_n) -
	f'(\sqrt{\frac{n-2}{n}} T_{n-1})] \leq A_1 + A_2 + A_3$ where
\begin{itemize}
\item $A_1 = \frac{ 3 \ee(|X_n|)}{\sqrt{{n \choose 2}}}$.

\item $A_2 = \frac{2 \sqrt{\frac{n-2}{n}}}{\sqrt{{n \choose 2}}}
\ee(|X_n||T_{n-1}|)$.

\item $A_3 = \ee \left[ \frac{|X_n|}{b \sqrt{{n \choose 2}}} \int_0^1
\ii_{[z,z+b]}(\sqrt{\frac{n-2}{n}} T_{n-1}+
\frac{s X_n}{\sqrt{{n \choose 2}}}) ds \right]$.
\end{itemize}  By part 3 of Lemma
	\ref{ismart}, $\ee(|X_n|) \leq \sqrt{\ee(X_n^2)} =
	\sqrt{n-1}$; thus $A_1 \leq \frac{3 \sqrt{2}}{\sqrt{n}}$. By
	parts 3 and 4 of Lemma \ref{ismart}, \[ \ee(|X_n| |T_{n-1}|)
	\leq \sqrt{\ee(X_n^2) \ee(T_{n-1}^2)} = \sqrt{n-1}. \] Thus
	$A_2 \leq \frac{2 \sqrt{2}}{\sqrt{n}}$.

	Note that $A_3=A_3'+A_3''$ where \[ A_3' = \ee \left[
	\ii(|X_n| \leq 2e \sqrt{n}) \frac{|X_n|}{b \sqrt{{n \choose
	2}}} \int_0^1 \ii_{[z,z+b]}(\sqrt{\frac{n-2}{n}} T_{n-1}+
	\frac{s X_n}{\sqrt{{n \choose 2}}} ) ds \right] \] and \[
	A_3'' = \ee \left[ \ii(|X_n| > 2e \sqrt{n}) \frac{|X_n|}{b
	\sqrt{{n \choose 2}}} \int_0^1
	\ii_{[z,z+b]}(\sqrt{\frac{n-2}{n}} T_{n-1}+ \frac{s
	X_n}{\sqrt{{n \choose 2}}}) ds \right].\] Clearly
	\[ A_3' \leq \frac{2e \sqrt{n}}{b \sqrt{{n \choose 2}}} \ee
	\left[ \ii_{[z-\frac{2e \sqrt{n}}{\sqrt{{n \choose
	2}}},z+b+\frac{2e \sqrt{n}}{\sqrt{{n \choose 2}}}]} \left(
	\sqrt{\frac{n-2}{n}} T_{n-1} \right) \right].\] Now use the
	fact (explained in \cite{Man}) that \[ 0 \leq \ee(\ii_B(c_1
	T_n + c_2)) \leq \frac{|B|}{c_1 \sqrt{2 \pi}} + 2 \delta(0,n)
	\] for any interval $B$ and constants $c_1,c_2$ with $c_1 \neq
	0$. It follows that $A_3' \leq \frac{D_1}{\sqrt{n}} +
	\frac{D_2}{bn} + \frac{D_3 \delta(0,n-1)}{b \sqrt{n}}$ where
	$D_1,D_2,D_3$ are universal constants. To bound $A_3''$, note
	that since $|X_n| \leq n$, one has that $A_3'' \leq \frac{n}{b
	\sqrt{{n \choose 2}}} \pp(|X_n| > 2e \sqrt{n})$. The proof
	of Proposition 4.6 of \cite{F1} derives the concentration
	inequality $\pp(|X_n| > 2e \sqrt{n}) \leq 2e^{-2e
	\sqrt{n}}$. Since $b$ will later be chosen to be a constant
	multiplied by $n^{-1/2}$, it follows that $A_3''$ is much
	smaller than $A_3'$ for large $n$, and one concludes that
	\[ A_3 \leq \frac{D_1}{\sqrt{n}} + \frac{D_2}{bn} + \frac{D_3
	\delta(0,n-1)}{b \sqrt{n}}\] where $D_1,D_2,D_3$ are universal
	constants.

	Combining the bounds on $A_1,A_2,A_3$, we conclude that \[ \ee
	\left[ f'(T_n) - f'(\sqrt{\frac{n-2}{2}} T_{n-1}) \right] \leq
	\frac{D_1}{\sqrt{n}} + \frac{D_2}{bn} + \frac{D_3
	\delta(0,n-1)}{b \sqrt{n}}\] where $D_1,D_2,D_3$ are universal
	constants.

	Next, we upper bound \[ \ee \left| \frac{ X_n^2}{n-1} \int_0^1
	\left[f'(\sqrt{\frac{n-2}{n}} T_{n-1} + t \frac{X_n}{\sqrt{{n
	\choose 2}}}) - f'(\sqrt{\frac{n-2}{n}} T_{n-1}) \right] dt \right|.\]
	Arguing as in the previous paragraph this is at most
	$B_1+B_2+B_3$ where
\begin{itemize}
\item $B_1 = \frac{1}{n-1} \int_{0}^1 \frac{3t \ee(|X_n|^3)}{\sqrt{{n
\choose 2}}} dt = \frac{3 \ee(|X_n|^3)}{2 (n-1) \sqrt{{n \choose 2}}}$.

\item $B_2 = \frac{1}{n-1} \int_{0}^1 \frac{2t
\sqrt{\frac{n-2}{n}}}{\sqrt{{n \choose 2}}} \ee(|T_{n-1}||X_n|^3) dt =
\frac{\sqrt{\frac{n-2}{n}}}{(n-1) \sqrt{{n \choose 2}}}
\ee(|T_{n-1}||X_n|^3)$.

\item $B_3 = \frac{1}{n-1} \ee \left[ \frac{|X_n|^3}{b \sqrt{n \choose 2}} \int_{0}^1 \int_0^1 t \ii_{[z,z+b]}(\sqrt{\frac{n-2}{n}} T_{n-1} + st \frac{X_n}{\sqrt{{n
\choose 2}}}) ds dt \right]$.
\end{itemize} To bound $B_1$, use part 1 of Lemma \ref{boundmom} to conclude that $B_1 \leq \frac{3}{\sqrt{n}}$. To
 bound $B_2$, use part 2 of Lemma \ref{boundmom} to conclude that $B_2
\leq \frac{2}{\sqrt{n}}$. To bound $B_3$, one uses an argument almost
identical to that for $A_3$ to conclude that  \[ B_3 \leq
\frac{E_1}{\sqrt{n}} + \frac{E_2}{bn} + \frac{E_3 \delta(0,n-1)}{b
\sqrt{n}}\] where $E_1,E_2,E_3$ are universal constants.

	Summarizing, it has been proved that \[ \delta(b,n) \leq
\frac{C_1}{\sqrt{n}} + \frac{C_2}{bn} + \frac{C_3 \delta(0,n-1)}{b
\sqrt{n}}\] where $C_1,C_2,C_3$ are universal constants. From
\cite{Bol} or \cite{Man}, $\delta(0,n) \leq \delta(b,n) +
\frac{b}{\sqrt{2 \pi}}$ for all $b$, which implies that \[ \delta(0,n)
\leq \frac{C_1}{\sqrt{n}} + \frac{C_2}{bn} + \frac{C_3
\delta(0,n-1)}{b \sqrt{n}} + \frac{b}{\sqrt{2 \pi}}.\] We argue by
induction that there is a universal constant $C$ so that $\delta(0,n)
\leq \frac{C}{\sqrt{n}}$ for all $n$. Assuming the result for $n-1$,
one obtains that
\[ \delta(0,n) \leq \frac{C_1}{\sqrt{n}} + \frac{C_2}{bn} + \sqrt
{\frac{3}{2}} \frac{C \cdot C_3 }{b n} + \frac{b}{\sqrt{2 \pi}}.\]
Choosing $b= \frac{2 C_3}{\sqrt{n}}$, it follows that if $C$ is
sufficiently large, the induction step will work for all $n$. This
completes the proof. \end{proof}

	To conclude this section, we note that it would be of interest
	to prove the following (more general) conjecture. An error
	term of $O(n^{-1/4})$ has recently been established by the
	method of exchangeable pairs \cite{F4}.

\begin{center}
\end{center}

{\bf Conjecture}: Let $i \geq 2$ be fixed. Then for all $n \geq i$ and real $x_0$, \[ \left| \pp \left( \sqrt{\frac{n!}{(n-i)! i}}
\frac{\chi^{\lambda}(12 \cdots i)}{dim(\lambda)} \leq x_0 \right) -
\frac{1}{\sqrt{2 \pi}} \int_{-\infty}^{x_0} e^{-\frac{t^2}{2}} dt \right|
\leq C_i n^{-1/2} \] where $C_i$ is a constant depending on $i$.

\section{Central limit theorem for Jack measure} \label{jack}

	For $\alpha>0$ the Jack$_{\alpha}$ measure on partitions of
size $n$ chooses a partition $\lambda$ with probability
\[\frac{\alpha^n n!}{\prod_{x \in \lambda} (\alpha a(x) + l(x) +1)
(\alpha a(x) + l(x) + \alpha)}, \] where the product is over all boxes
in the partition. Here $a(x)$ denotes the number of boxes in the same
row of $x$ and to the right of $x$ (the ``arm'' of x) and $l(x)$
denotes the number of boxes in the same column of $x$ and below $x$
(the ``leg'' of x). For example the partition of 5 below \[
\begin{array}{c c c} \framebox{\ }& \framebox{\ }& \framebox{\ } \\
\framebox{\ }& \framebox{\ }& \end{array} \] would have
Jack$_{\alpha}$ measure \[\frac{30 \alpha^2}{(3 \alpha+1) (\alpha+2)(2
\alpha+1)(\alpha+1)^2}.\] Note that when $\alpha=1$, Jack measure
reduces to Plancherel measure of the symmetric group. The papers
\cite{O2}, \cite{BO1} emphasize that for $\alpha$ fixed the study of
Jack$_{\alpha}$ measure is an important open problem, about which
relatively little is known for general values of $\alpha$. It is a
discrete analog of eigenvalue ensembles from random matrix theory and like Jack
polynomials \cite{GHJ}, should also be relevant to the moduli space of
curves.

	Given $\alpha>0$, the quantity to be studied is \[
	T_{n,\alpha}(\lambda) = \frac{\sum_i (\alpha {\lambda_i
	\choose 2} - {\lambda_i' \choose 2})}{\sqrt{ \alpha {n \choose
	2}}},\] where as usual $\lambda_i$ is the length of the ith
	row of $\lambda$ and $\lambda_i'$ is the length of the ith
	column of $\lambda$. It is of interest to study the quantity
	$T_{n,\alpha}(\lambda)$ under Jack measure for several reasons. When
	$\alpha=1$ it reduces to the study of the character ratio of
	transpositions under Plancherel measure. When $\alpha=2$ it is
	a spherical function of the Gelfand pair $(S_{2n},H_{2n})$
	where $H_{2n}$ is the
	hyperoctahedral group of size $2^nn!$. Also by Corollary 1 of
	\cite{DHol}, there is a natural random walk on perfect
	matchings of the complete graph on $n$ vertices, whose eigenvalues are
	precisely $\frac{T_{n,2}(\lambda)}{\sqrt{n(n-1)}}$, occurring
	with multiplicity proportional to the Jack$_2$ measure of
	$\lambda$.

	The paper \cite{F2} used the ``exchangeable pairs'' version of
	Stein's method to prove a central limit theorem for
	$T_{n,\alpha}$ with error term $C_{\alpha} n^{-1/4}$ where
	$C_{\alpha}$ is a constant depending on $\alpha$. This was
	sharpened in \cite{F3} using martingales to $C_{\alpha,s}
	n^{-s}$ for any $s<\frac{1}{2}$.

	The main result of this section is Theorem \ref{main2}.

\begin{theorem} \label{main2} Suppose that $\alpha \geq 1$ and let $\lambda$ be chosen from the Jack$_{\alpha}$ measure on partitions of size $n$. Then there is a constant $C_{\alpha}$
	depending on $\alpha$ so that for all $n \geq 2$ and real
	$x_0$, \[ \left| \pp \left( T_{n,\alpha}(\lambda) \leq x_0
	\right) - \frac{1}{\sqrt{2 \pi}} \int_{-\infty}^{x_0} e^{-
	\frac{t^2}{2}} dt \right| \leq C_{\alpha} n^{-1/2}. \]
	\end{theorem}

	Note that in Theorem \ref{main2} we suppose that $\alpha \geq
	1$ since the Jack$_{\alpha}$ probability of $\lambda$ is equal
	to the Jack$_{\frac{1}{\alpha}}$ probability of the transpose
	of $\lambda$, implying that for any $x$, the Jack$_{\alpha}$
	probability that $T_{n,\alpha}=x$ is equal to the
	Jack$_{\frac{1}{\alpha}}$ probability that
	$T_{n,\frac{1}{\alpha}}=-x$. Also note that $C_{\alpha}$ must
	depend on $\alpha$, since by Corollary 5.3 of \cite{F2}, the
	random variable $T_{n,\alpha}$ has mean 0, variance 1, and
	third moment $\frac{\alpha-1}{\sqrt{\alpha {n \choose 2}}}$.

	There is no need to write out a proof of Theorem \ref{main2},
	which uses exactly the same logic as that of Theorem
	\ref{main1}. But it is necessary to give analogs of Lemmas
	\ref{ismart}, \ref{moments}, \ref{identity}, and
	\ref{boundmom}, and we do that. The concentration inequality
	needed for $X_{n,\alpha}$ is Lemma 6.6 of \cite{F2} and can be
	used to give another proof of some results in this section.

	There is an $\alpha$-analog of Kerov's growth process (due to
Kerov \cite{K4}) giving a sequence of partitions
$(\lambda(1),\cdots,\lambda(n))$ with $\lambda(j)$ distributed
according to the Jack$_{\alpha}$ measure on partitions of size $j$;
see \cite{F3} for details. Moreover from the definition of
$T_{n,\alpha}$, it follows that \[ T_{n,\alpha} =
\frac{1}{\sqrt{\alpha {n \choose 2}}} (X_{1,\alpha} + \cdots +
X_{n,\alpha}). \] Here $X_{1,\alpha}=0$ and if $j \geq 2$ then
$X_{j,\alpha}=c_{\alpha}(x)$ where $x$ is the box added to
$\lambda(j-1)$ to obtain $\lambda(j)$ and the ``$\alpha$-content''
$c_{\alpha}(x)$ of a box $x$ is defined to be $\alpha$ (column number
of x-1) - (row number of x-1).

	Lemma \ref{ismart2} is an analog of Lemma \ref{ismart} and
	is generalized in \cite{F3} to arbitrary spherical functions of the
	Gelfand pair $(S_{2n},H_{2n})$.

\begin{lemma} \label{ismart2} (\cite{F3})
\begin{enumerate}
\item $\ee(X_{j,\alpha}|\lambda(j-1))=0$ for $2 \leq j \leq n$ and all partitions $\lambda(j-1)$.
\item $\ee(X_{j,\alpha}|T_{n,\alpha})=\frac{(j-1) \sqrt{\alpha}}{\sqrt{ {n \choose 2}}} T_{n,\alpha}$ for all $1 \leq j \leq n$.
\item $\ee(X_{j,\alpha}^2) = \alpha(j-1)$.
\item $\ee(T_{n,\alpha}^2) = 1$.
\end{enumerate}
\end{lemma}

	Lemma \ref{moments2} is the $\alpha$ version of Lemma \ref{moments}.

\begin{lemma} \label{moments2} Let $\lambda(j-1)$ be a partition of size $j-1 \geq 1$.
\begin{enumerate}
\item (\cite{K4}) $\ee(X_{j,\alpha}^2|\lambda(j-1))= \alpha(j-1)$.
\item (\cite{La}) \begin{eqnarray*} \ee(X_{j,\alpha}^4|\lambda(j-1)) & = & \alpha^2 {j \choose 2} + \alpha(\alpha-1)^2(j-1) + 3 \alpha \sum_{x \in \lambda(j-1)} c_{\alpha}(x)^2 \\  & & + 3\alpha(\alpha-1) \sum_{x \in \lambda(j-1)} c_{\alpha}(x). \end{eqnarray*}
\end{enumerate}
\end{lemma}

	Lemma \ref{identity2} is the $\alpha$ version of Lemma
	\ref{identity}. The proof is combinatorial, as opposed to the
	algebraic argument given for Lemma \ref{identity}.

\begin{lemma} \label{identity2} Consider the Jack$_{\alpha}$ measure on partitions of size $n$.
\begin{enumerate}
\item If $m \geq 1$ is an integer then \[ \ee \left( \prod_{x \in \lambda} (m+c_{\alpha}(x)) \right) = m^n.\]
\item Let $e_{r,\alpha}(\lambda)$ denote the $r$th elementary symmetric
function of the $\alpha$-contents of the boxes of $\lambda$. Then $\ee(e_{r,\alpha}(\lambda))=0$ for $1 \leq r \leq n$.
\end{enumerate}
\end{lemma}

\begin{proof} It suffices
 to prove the first assertion since the second assertion follows from
 the first by taking the coefficient of $m^{n-r}$ on both sides. Page
 324 of \cite{Mac} proves the identity \[ \sum_{\lambda}
 b_{\lambda}(q,t) P_{\lambda}(y;q,t) P_{\lambda}(z;q,t) = e^{ \sum_{n
 \geq 1} ( \frac{1}{n} \frac{1-t^n}{1-q^n} p_n(y) p_n(z))}\] where the
 sum is over all $\lambda$ of all sizes, $P_{\lambda}(y;q,t)$ denotes
 a Macdonald symmetric function, $p_n(y) = \sum_i y_i^n$ denotes the
 nth power sum symmetric function, and $b_{\lambda}(q,t)$ is a number
 to be discussed more below. We apply the homomorphism of the ring of
 symmetric functions determined by $p_r(y) \mapsto mu^r$, $p_r(z)
 \mapsto l^{1-r}$ for all $r \geq 1$ where $m,l$ are positive
 integers; this is possible since the $p_r$'s are algebraically
 independent. Then we take the limit $q=t^{\alpha},t \mapsto 1$ in
 which Macdonald polynomials become Jack polynomials.

With these substitutions, consider the left hand side of the
        identity. By pages 380 and 381 of \cite{Mac}, \[
        b_{\lambda}(q,t) \mapsto \prod_{x \in \lambda} \frac{\alpha
        a(x) + l(x)+1}{\alpha a(x) + l(x) + \alpha} \]
        \[P_{\lambda}(y;q,t) \mapsto u^{|\lambda|} \prod_{x \in
        \lambda} \frac{m+c_{\alpha}(x)}{\alpha a(x) + l(x)+1} \]
        \[ P_{\lambda}(z;q,t) \mapsto \frac{1}{l^{|\lambda|}} \prod_{x
        \in \lambda} \frac{l+c_{\alpha}(x)}{\alpha a(x) + l(x)
        +1}. \] Letting $l \rightarrow \infty$, one sees that the
        coefficient of $u^n$ in the left-hand side of the identity is
        $\frac{1}{\alpha^n n!} \ee \left( \prod_{x \in \lambda}
        (m+c_{\alpha}(x)) \right)$.

	 Consider the right hand side of the identity with these
        substitutions. One obtains $e^{\sum_{n \geq 1}
        (\frac{mu^n}{n \alpha l^{n-1}})}$. Letting $l \rightarrow \infty$,
        one obtains $e^{\frac{mu}{\alpha}}$, and taking the
        coefficient of $u^n$ gives $\frac{m^n}{\alpha^n
        n!}$. Comparing with the previous paragraph proves the first
        assertion of the lemma. \end{proof}

	Finally, we give the analog of Lemma \ref{boundmom}.	

\begin{lemma} \label{boundmom2} Suppose that $n \geq 3$. There is a constant $D_{\alpha}$ such that
\begin{enumerate}
\item $\ee(|X_{n,\alpha}|^3) \leq D_{\alpha} n^{3/2}$.
\item $\ee(|T_{n-1,\alpha}||X_{n,\alpha}|^3) \leq D_{\alpha} n^{3/2}$.
\end{enumerate} \end{lemma}

\begin{proof} The proof
 method is the same as that of Lemma \ref{boundmom}, using the
 Cauchy-Schwarz inequality in the first part and the conditional
 Cauchy-Schwarz inequality in the second part. One uses that
 $\ee (X_{n,\alpha}^2|\lambda(n-1))=\alpha(n-1)$ for all
 $\lambda(n-1)$ (part 1 of Lemma \ref{moments2}). Also one needs that
\begin{eqnarray*}
\ee(X_{n,\alpha}^4) & = & \alpha^2 {n \choose 2} +
\alpha(\alpha-1)^2(n-1) + 3 \alpha \ee \left( \sum_{x \in
\lambda(n-1)} c_{\alpha}(x)^2 \right)\\ & = & \alpha^2 {n \choose 2} +
3 \alpha^2 {n-1 \choose 2} + \alpha(\alpha-1)^2(n-1). \end{eqnarray*}
The first equality used part 2 of Lemma \ref{moments2} and the fact
that $\ee(T_{n-1,\alpha}) = 0$. The second equality used part 4 of Lemma
\ref{ismart2} and Lemma \ref{identity2} with $r=2$. \end{proof}

\section{Acknowledgements} The author was partially supported by NSA
	grant number H98230-05-1-0031.


\begin{thebibliography}{AAA}

\bibitem [AlD] {AD} D. Aldous and P. Diaconis, Longest increasing
subsequences: from patience sorting to the Baik-Deift-Johansson
theorem, Bull. AMS (N.S.) {\bf 36} (1999) 413-432.

\bibitem [Bol]{Bol} E. Bolthausen, An estimate of the remainder term
in a combinatorial central limit theorem, {Z. Wahrsch. Verw. Gebiete}
{\bf 66} (1984) 379-386.

\bibitem [BOO]{BOO} A. Borodin, A. Okounkov, and G. Olshanski,
Asymptotics of Plancherel measures for symmetric groups,
J. Amer. Math. Soc. {\bf 13} (2000) 481-515.

\bibitem [BO] {BO1} A. Borodin and G. Olshanski, Z-measures on
partitions and their scaling limits, European J. Combin. {\bf 26}
(2005) 795-834.

\bibitem [De]{De} P. Deift, Integrable systems and combinatorial
theory, Notices Amer. Math. Soc. {\bf 47} (2000) 631-640.

\bibitem [DG]{DG} P. Diaconis and C. Greene, Applications of
Murphy's elements, Stanford University technical report no. 335
(1989).

\bibitem [DHol]{DHol} P. Diaconis and S. Holmes, Random walk on
trees and matchings, Elec. J. Probab. {\bf 7} (2002), 17 pages
(electronic).

\bibitem [DSh]{DS} P. Diaconis and M. Shahshahani, Generating a
random permutation with random transpositions,
Z. Wahr. Verw. Gebiete {\bf 57} (1981) 159-179.

\bibitem [EO]{EO} A. Eskin and A. Okounkov, Asymptotics of branched
coverings of a torus and volumes of moduli spaces of holomorphic
differentials, Invent. Math. {\bf 145} (2001) 59-103.

\bibitem [Fr]{Fr} F. Frobenuis, Uber die charaktere der symmetrischen
gruppe, Sitz. Konig. Preuss. Akad. Wissen. (1900) 516-534;
Gesammelte abhandlungen III, Springer-Verlag, Heidelberg, 1968,
148-166.

\bibitem [F1]{F1} J. Fulman, Stein's method and Plancherel measure of
the symmetric group, Trans. Amer. Math. Soc. {\bf 357} (2005)
555-570.

\bibitem [F2]{F2} J. Fulman, Stein's method, Jack measure, and the
Metropolis algorithm, J. Combin. Theory Ser. A {\bf 108} (2004)
275-296.

\bibitem [F3]{F3} J. Fulman, Martingales and character ratios,
Trans. Amer. Math. Soc. {\bf 358} (2006) 4533-4552.

\bibitem [F4]{F4} J. Fulman, Stein's method and random character
ratios, to appear in Trans. Amer. Math. Soc., available at
http://www-rcf.usc.edu/$\sim$fulman.

\bibitem [GHJ]{GHJ} I. Goulden, J. Harer, and D. Jackson, A geometric
parametrization for the virtual Euler characteristic of the moduli
spaces of real and complex algebraic curves,
Trans. Amer. Math. Soc. {\bf 353} (2001) 4405-4427.

\bibitem [Ho]{Ho} A. Hora, Central limit theorem for the adjacency
operators on the infinite symmetric group, Comm. Math. Phys.
{\bf 195} (1998) 405-416.

\bibitem [IO]{IO} V. Ivanov and G. Olshanski, Kerov's central limit
theorem for the Plancherel measure on Young diagrams, in: Symmetric
Functions 2001: Surveys of developments and perspectives, Kluwer
Academic Publishers, S. Fomin, Editor, Dodrecht, 2002, pp. 93-151.

\bibitem [J]{J} K. Johansson, Discrete orthogonal polynomial
ensembles and the Plancherel measure, Ann. of Math. (2) {\bf
153} (2001) 259-296.

\bibitem [K1]{K1} S.V. Kerov, Gaussian limit for the Plancherel
measure of the symmetric group, Compt. Rend. Acad. Sci. Paris,
Serie I, {\bf 316} (1993) 303-308.

\bibitem [K2]{K2} S.V. Kerov, The boundary of Young lattice and
random Young tableaux, in:  Formal power series and algebraic
combinatorics, DIMACS Ser. Discrete Math. Theoret. Comput. Sci. {\bf
24}, Amer. Math. Soc., Providence, RI, 1996, pp. 133-158.

\bibitem [K3]{K3} S.V. Kerov, Transition probabilities of continual
Young diagrams and the Markov moment problem, Funct. Anal. Appl.
{\bf 27} (1993) 104-117.

\bibitem [K4]{K4} S.V. Kerov, Anisotropic Young diagrams and Jack
symmetric functions,  Funct. Anal. Appl. {\bf 34} (2000) 41-51.

\bibitem [La]{La} M. Lassalle, Jack polynomials and some identities
for partitions, Trans. Amer. Math. Soc. {\bf 356} (2004)
3455-3476.

\bibitem [Mac]{Mac} I. Macdonald, Symmetric functions and Hall
polynomials, Second edition, Oxford University Press, New York, 1995.

\bibitem [Man]{Man} B. Mann, Bolthausen's proof of Berry-Esseen,
unpublished manuscript (1994).

\bibitem [Mu]{Mu} G.E. Murphy, A new construction of Young's
seminormal representation of the symmetric group, J. Algebra
{\bf 69} (1981) 287-291.

\bibitem [O1]{O1} A. Okounkov, Random matrices and random
permutations, Internat. Math. Res. Notices {\bf 20} (2000)
1043-1095.

\bibitem [O2]{O2} A. Okounkov, The uses of random partitions, in: XIV
International Congress on Mathematical Physics, World Sci. Publ., 2005,
pp. 379-403.

\bibitem [OP]{OP} A. Okounkov and R. Pandaripandhe, Gromov-Witten
theory, Hurwitz numbers, and Matrix models, I, arXive:math.AG/0101147.

\bibitem[Sa]{Sa} B. Sagan, The symmetric group. Representations,
combinatorial algorithms, and symmetric functions, Springer-Verlag,
New York, 1991.

\bibitem[ShSu]{SS} Q. Shao and Z. Su, The Berry-Esseen bound for
character ratios, Proc. Amer. Math. Soc. {\bf 134} (2006), 2153-2159.

\bibitem[Sn1]{Sn1} P. Sniady, Asymptotics of characters of symmetric
groups, Gaussian fluctuations of Young diagrams and genus expansion,
arXive:math.CO/0411647.

\bibitem[Sn2]{Sn2} P. Sniady, Gaussian fluctuations of characters of
symmetric groups and of Young diagrams, arXive:math.CO/0501112.

\bibitem[Stn]{Stn} C. Stein, Approximate computation of expectations,
Institute of Mathematical Statistics Lecture Notes, Volume 7, 1986.

\end{thebibliography}
 \end{document}